\documentclass[12pt]{article}
 \usepackage{amsfonts,amssymb,eucal}
\usepackage{amsthm}
 \usepackage{amsmath}
 \usepackage{amscd}
 \usepackage{latexsym}

\usepackage{epsfig}
 \usepackage{rotating}
 \usepackage{psboxit}

 \numberwithin{equation}{section}

 \newcommand{\got}[1]{{\mathfrak{#1}}}% gothic with mbox for  mathematic
 \newcommand{\gata}{\blacksquare}

\def\openone{\leavevmode\hbox{\small1\kern-3.8pt\normalsize1}}%
                  \newcommand{\oo}{\openone}
\newcommand{\R}{\ensuremath{\mathbb{R}}}
\newcommand{\C}{\ensuremath{\mathbb{C}}}

 \renewcommand{\P}{\ensuremath{\mathbb{P}}}
 
 \newcommand{\Hi}{\ensuremath{\mathcal{H}}}% Hilbert space
 \renewcommand{\P}{\ensuremath{\mathbb{P}}}
 \newcommand{\Gras}{\mbox{$G_n({\C}^{m+n})$}}% Grassmannian
 
 \newtheorem{theorem}{Theorem}
\newtheorem{lemma}{Lemma}
 \newtheorem{Remark}{Remark}
\newtheorem{proposition}{Proposition}
\newcommand{\gcm}{\ensuremath{\mathbb{C}}}% C double for exponent
\newcommand{\skp}[2]{{\langle #1,#2\rangle}}
\newcommand{\df}{\mbox{$:=$}}% def
\begin{document}

\vspace*{2.5cm}
\noindent

 \begin{center} 
{\Large{\bf{ Geometrical  phases on 
hermitian symmetric spaces }}}\\[3ex]
Stefan Berceanu\\[2ex]
 National 
 Institute for Physics and Nuclear Engineering\\
         Department of Theoretical Physics\\
         PO BOX MG-6, Bucharest-Magurele, Romania\\
         E-mail: Berceanu@theor1.theory.nipne.ro\\
 \end{center}

\begin{abstract}
\footnotesize For simple Lie groups, 
  the only   homogeneous   manifolds $G/K$, where
$K$ is maximal compact subgroup,
 for which the phase of the scalar product of two
coherent state vectors is twice the symplectic area 
of a geodesic triangle are the hermitian symmetric spaces.
An explicit calculation of the multiplicative factor
 on the complex Grassmann manifold and its noncompact dual is presented.
 It is shown that the multiplicative factor  is identical with the
 two-cocycle  considered by A. Guichardet and
D.  Wigner for simple Lie groups. 
\end{abstract}
\normalsize

\section{Introduction}
   Six questions referring to the relationship between  coherent 
states and geometry have been  presented in \cite{clasic}. In the same context,
 reference  \cite{sb7} 
was devoted to the following 
question: find a geometric significance of the phase of the
 scalar product of two
 coherent states. An explicit answer to this question  for the Riemann sphere
 was given by Perelomov  (cf. reference \cite{per1}). Earlier,
S.  Pancharatnam \cite{pan,berry} showed that the phase difference between
 the
initial and final state is $<A|A'>=\exp (-i\Omega_{ABC}/2)$, where
$\Omega_{ABC}$ is the solid angle subtended by the geodesic triangle $ABC$
on the Poincar\'e sphere.
 The holonomy of a loop in the projective Hilbert space is twice the
symplectic area of any two-dimensional submanifold whose boundary is the given
loop (see  \cite{aa} and  Proposition 5.1 in 
\cite{ma}).

A general answer to the question of the geometric significance of
 the phase of the scalar product of two coherent state vectors using
the coherent state embedding and the so called ``Cauchy formulas''
 was given in \cite{sb7} and \cite{martin}. In reference \cite{sb7}
it was proved that for compact hermitian symmetric spaces the phase 
 of the scalar product of two coherent states is twice the symplectic 
area of a geodesic triangle determined by the corresponding points
 on the manifold and the origin of the system of coordinates. A
similar 
result was also obtained in another formulation in \cite{clerc}.
In fact, in reference \cite{sb7} this result was proved on a restricted class
of manifolds: the compact, homogeneous, simply connected Hodge manifolds, which
are in the same time naturally reductive.
But  this class of manifolds
considered by me  in \cite{sb7} consists in
fact only of  the Hermitian symmetric spaces \cite{martinb}. Indeed,
any naturally reductive space with an invariant almost K\"ahler
structure is  locally  Hermitian
symmetric (cf. Corollary 7 in \cite{van}; see also
Corollary 9 in the same reference) and simply connectedness implies
 Hermitian symmetry.
 On the other side, the results of reference \cite{sb7}
  are also true for  other
manifolds than compact Hermitian symmetric spaces.
 For example, the results are true for
 the Heisenberg-Weyl group \cite{per1} as well as for the
noncompact dual of the  complex Grassmann manifold \cite{brasov}.    
 
An explicit
formula was presented for the symplectic area of geodesic triangles on the
complex Grassmann manifold \cite{sb7},
 and also for its noncompact dual (\cite{brasov},
\cite{montreal}).  Lately, I learned that  the formula for the symplectic
 area on the noncompact Grassmann manifold was  found out earlier in the
 paper 
\cite{domic} devoted to the Gromov's norms.
 The methods of reference \cite{domic} were developed in
\cite{clerc}. Also there are other  references  on two-cocycles
on real simple Lie groups, which are related to the symplectic
area of geodesic triangles \cite{guit,dupont,guit1}. 

A. Guichardet and D. Wigner \cite{guit}
 have  proved that a simple Lie group has non-trivial
 continuous 2-cohomology group $H^2(G, \R )$ if and only if $G/K$ admits a
 $G$-invariant complex structure,  where $K$ is
  a maximal compact subgroup of $G$.
In this context, let us remained some well known facts 
(cf. e.g. \cite{cal}):
If $\got{g}$ is the Lie algebra of the compact and connected Lie group $G$,
then $H^q(\got{g})$ is isomorphic with the $q^{\text{th}}$
 cohomology group $H^q(G)$ with
 real 
coefficients and the ring $H(\got{g})$ is isomorphic with the cohomology
ring $H(G)$ of $G$.
 If $\got{g}$ is a semi-simple Lie algebra over a field of characteristic 0,
 then
 $H^1 (\got{g})=\{ 0 \}$, $H^2 (\got{g})=\{ 0 \}$ and $H^3(\got{g})\not=
 \{ 0\}$.  Moreover,  for a simply connected Lie group $G$, not 
only $H^1(G)=\{ 0 \}$, but  also $H^2(G)=\{ 0 \}$.

In the present paper we give an explicit calculation of 
 the multiplicative factor of
 representations on the
 complex Grassmann manifold, which, when expressed in
Pontrjagin's coordinates, it  is shown to be identical with the two-cocycle
considered by A. Guichardet and D. Wigner.
 The notation and technique for manipulating
the Grassmann manifold  are that from references
\cite{sbl}, \cite{grass}.

The geometric significance of the 2-cocycle (see below eq.
(\ref{2coc})) as a symplectic
 area of a geodesic triangle was found by
J.-L.  Dupont and A.  Guichardet \cite{dupont}. Using the results of 
\cite{guit,dupont,guit1} and our results in \cite{sb7,brasov}, it
follows that: {\em If $G $  is  a simple Lie group and $K$ a maximal compact
 subgroup, then the only  coherent state manifolds
$G/K$ for which the phase of the scalar product of two
coherent state vectors is twice the symplectic area 
of a geodesic triangle are the hermitian symmetric spaces.} This remark is
a completion of our assertions in \cite{sb7}. 

In this context, the following  question naturally arise:
 {\em For which Lie groups $G$,  which admits coherent state
representations} (cf. \cite{lis1}, \cite{neeb}), {\em the assertion   
``the phase of the scalar product
of two coherent states is twice the symplectic area of geodesic triangles''
 is still true?}

Let us also remained that 
G. Lion and M. Vergne \cite{lion} have underlined the
 relationship between the 2-cocycle of the Segal-Shale-Weil projective 
 representation of
 the symplectic group $G=Sp(B)$ of the vector space  $(V,B)$ and the Maslov
 index. In the same context we mention also the work of
B. Magneron  \cite{magn}.

 The problem of symplectic
 area of geodesic triangles on symmetric spaces
was considered  also by A. Weinstein \cite{wein}. See also
the paper \cite{pierre}. Let us mention also the paper \cite{cefacem}.

The present paper is laid out as follows: in \S \ref{br}
simplest  manifolds on which the phase of the scalar product is twice the
 symplectic are  presented
--- the sphere $SU(2)/U(1)$, its noncompact dual $SU(1,1)/U(1)$, and
the Heisenberg group --- while in \S \ref{question}
 are
reviewed our own results referring to
 the Grassmann manifold. In \S \ref{two} we present
a calculation of 
   the phase which appears when we multiply
two representations on the
 complex Grassmann manifold and its noncompact dual.  Our 
results are compared with those in references \cite{guit},
\cite{dupont},\cite{guit1} in \S \ref{compar}. The necessary
  formulas referring
to the complex Grassmann manifold and its noncompact dual
 are collected  in \S \ref{app}. The definitions of coherent states are
as in references \cite{tim,buc}. 

\section{Previous results}\label{prev}
\subsection{The phase of the scalar product of coherent states =\\ (2$\times$)
symplectic area - basic examples}\label{br}

a) Let us consider  the sphere $S^2=SU(2)/U(1)=\C\P^1 $.
The commutation relations of the generators of the group $SU(2)$ are 
\begin{equation}\label{com1}
 [J_0,J_{\pm}]=\pm J_{\pm};~[J_-,J_+]=-2J_0. 
\end{equation}
We denote below with the same latter $X$ the generator of the Lie
algebra $\got{g}$ of the Lie group $G$
and the derived representation $d\pi (X)$ of the unitary irreducible
representation $\pi$ of the group $G$.
The action of the generators on the minimal weight vector
 $e_0$ ($e_0=|j,-j>$) is
$$ J_+e_0\not= 0; J_0e_0=-je_0; J_-e_0=0.  $$
The coherent state vectors are
$$e_z=e^{zJ_+}e_0,~z\in\C,$$
and the scalar product is
$$ (e_{\bar{z}},e_{\bar{z}'})=(1+z\bar{z}')^{2j}
=|~\cdot ~|e^{i\phi},$$
where the phase $\phi$ is
\begin{equation}\label{phi}
\phi = \frac{j}{2i}\log \frac{(1+z\bar{z}')}{(1+z'\bar{z})}~.
\end{equation}
It can be checked (caution: not an easy exercise!) that 
\begin{equation}\label{basic}
\boxed{\begin{array}{c}\text{the phase}~ \phi
\text{ of the scalar product of two coherent states}
 = \\(2\times) ~ {\mbox{\rm{symplectic area of the geodesic triangle}}}
\end{array}}
\end{equation}
where the two-form $\omega$ on the sphere is
$$  \omega = \frac{\sqrt{-1}}{2}
\frac{dz\wedge d\bar{z}}{(1+|z|^2)^2}~. $$
The    formula  (\ref{phi}) can be find in the book
of Perelomov  \cite{per1} at p. 63. See also 
 Pancharatnam \cite{pan}. 

%{\bf Exercise:} Check  that eq. (\ref{phi}) gives the spherical
% exceeding $\epsilon = 4\pi -\alpha 
%-\beta -\gamma$ in stereographic coordinates.

A similar formula with (\ref{phi}),
\begin{equation}\label{phi1}
\phi = -\frac{j}{2i}\log \frac{(1-z\bar{z}')}{(1-z'\bar{z})}~,
\end{equation}
holds for the noncompact manifold $SU(1,1)/U(1)$, with the same 
significance  (\ref{basic}). Here instead of (\ref{com1}) we have
the commutation relations
\begin{equation}\label{com2}
[K_0,K_{\pm}]=\pm K_{\pm}; [K_-,K_+]=2K_0.
\end{equation}
Taking  $e_0=|k,k>$, then
$$K_+e_0\not= 0,K_0e_0=ke_0, K_-e_0=0, k=1,\frac{3}{2},2,\frac{5}{2},\dots $$
The coherent states are
$$e_z= e^{zK_+}e_0; (e_{\bar{z}},e_{\bar{z}'})=(1-z\bar{z}')^{-2k},$$
and the two-form $\omega$ is 
$$  \omega = \frac{\sqrt{-1}}{2}
\frac{dz\wedge d\bar{z}}{(1-|z|^2)^2}~. $$
b) Now we consider the Heisenberg-Weyl group $\approx \C$.
The canonical commutation relations are
$$[a,a^+]=1; a^+e_0\not= 0; a e_0 =0.$$
The Glauber's coherent states are
$$e_z=e^{za^+}e_0,$$
with the scalar product
$$\ (e_{\bar{z}},e_{\bar{z}'})=e^{z\bar{z}'}.$$
Then
$\Im (z\bar{z}')= 2\times~ {\mbox{\rm{area of 
the geodesic triangle}}} (0,z, z+z')$   (cf. eq. (1.1.17) p.10 in \cite{per1}).
%\newpage
\subsection{Symplectic area of geodesic triangles  on the complex\\
 Grassmann manifold
and its noncompact dual}\label{question}

We remember the formulas which generalize expression
(\ref{phi})  ((\ref{phi1}))
of the symplectic area
on the sphere in stereographic coordinates ($SU(1,1)/U(1)$)
  to the complex 
Grassmann manifold (respectively, its  noncompact dual).
 The notation
referring to the complex Grassmann manifold and its noncompact dual
is collected  in the Appendix  (cf. \cite{sbl}, \cite{grass}).

Let us denote the   compact Grassmann manifold \Gras ~
of $n$-planes in $\C^{n+m}$ by
 \begin{equation}
X_c=G_c/K=SU(n+m)/S(U(n)\times U(m))~, \label{xc}
\end{equation}
and its noncompact dual by
\begin{equation}
X_n=G_n/K=SU(n,m)/ S(U(n)\times U(m))~.\label{xn}
\end{equation}

The following two theorems are extracted from \cite{sb7,montreal,brasov}.

 \begin{theorem}\label{gr}{Let $z,z'\in{\cal V}_0\subset\Gras$
(resp. its noncompact dual (\ref{xn})) be described by the
Pontrjagin's coordinates $Z, Z'$. Let $\gamma (0,z,z')$ be the geodesic
triangle obtained by
joining $0, z, z'$. Then the symplectic area of the surface 
$\sigma (0,z,z')$
 of the geodesic triangle $\gamma (0,z,z')$ is given by
 \begin{equation}\label{ar}
 {\got I}(0,Z,Z')=\int_{\sigma (0,z,z')}\omega
 =\frac{\epsilon }{4i}\log\,\frac{\det(\oo+\epsilon Z{Z'}^+)}
 {\det(\oo+\epsilon {Z'}Z^+)}~.
 \end{equation}}
 \end{theorem}
$\epsilon =1$ ($\epsilon = -1 $) corresponds to the compact
(noncompact) manifold $X_c$ (respectively, $X_n$).

The main ingredients for proving theorem \ref{gr} in the compact case are
presented in \cite{sb7}. A similar calculation can be done in the noncompact 
case (cf.  \cite{brasov}). Here we just remember that the two-form
$\omega$ in eq. (\ref{ar}) is
\begin{equation}
\omega = \frac{i}{2}\text{Tr}[dZ(\oo_n +\epsilon Z^+Z)^{-1}\wedge dZ^+(\oo_m 
+\epsilon Z Z^+)^{-1}].
\end{equation}

Let us mention that in the
 context of Gromov's norm of the K\"ahler class of symmetric domains,
 A.  Domic and D. Toledo \cite{domic} have calculated the symplectic area of
 the geodesic
 triangle in the case of the noncompact dual of the Grassmann
 manifold. Their calculation is based on  the Stokes's  formula 
 \begin{equation}\label{toledo}
 \int_{\Delta}\omega = \int_{\gamma (Q,R)}d^{\gcm}\rho_P,
 \end{equation}
Here
$\Delta$ is a geodesic simplex in the bounded symmetric domain,
  which  has the vertices $P, Q, R$  and $\rho_P$ is the
 (unique) potential for the Bergmann metric, i.e. a function such that
 $dd^{\gcm} \rho _P= \omega $, with the following properties: A) $\rho_P(P)=0$;
 B) $\rho_P$ is invariant under the isotropy group of $P$; C) 
 $d^{\gcm}\rho_P=0$ on geodesics through $P$. We only want to stress that
\begin{Remark} The proof in \cite{sb7} of
eq. (\ref{ar}) using Stokes's formula
is equivalent with the
 calculation of A. Domic and D. Toledo \cite{domic}.
\end{Remark}

 Recall the definition of Perelomov's coherent state vectors:
\begin{equation}
{e}_{Z,j}=\exp\sum_{{\varphi}\in\Delta^+_n}(Z_{\varphi}F^+_{\varphi})
j ,~~~~\underline{{e}}_{Z,j}=({e}_{Z,j} ,
{e}_{Z,j})^{-1/2}{e}_{Z,j} ,
\label{z}
\end{equation}
\begin{equation}
{e}_{B,j}=\exp\sum_{{\varphi}\in\Delta^+_n}
(B_{\varphi}F^+_{\varphi}-{\bar
B}_{\varphi}F^-_{\varphi}) j , \label{b}
~~~{e}_{B,j}\df  \underline{{e}}_{Z,j}.
\end{equation}
 where $\Delta^+_n$ are the positive noncompact roots, $ Z\df (Z_ \varphi )
\in {\C}^d$  ($d$ = complex dimension of $M$) are the  local
coordinates in the maximal neighborhood
${\cal V}_0 \subset  M $, $F^+_{\varphi} j\neq  0,F^-_{\varphi} j = 0,
~  \varphi\in\Delta^+_n $, and $j$ is the extreme (here minimal, see
eq.
(\ref{pondere}) below)
weight vector of the 
representation. Note that $\mathcal{V}_0\approx X_n$ for the noncompact case.

\begin{theorem}\label{sb2}{Let $M$ be a  hermitian
 symmetric manifold.
 Let us consider on the
manifold of coherent states $M$ the Perelomov's coherent vectors 
(\ref{z}) in a
local chart, corresponding to the fundamental representation $\pi$.
Let us consider the points $Z, Z'\in {\cal V}_0\subset M $ such that
$0, Z, Z'$ is a geodesic triangle.
 Then the phase $\Phi_M$ defined by the relation
\begin{equation}\label{pufi}
(\underline{{e}}_{Z'},\underline{{e}}_{Z})=
|(\underline{{e}}_{Z'},\underline{{e}}_{Z})|\exp({i\Phi_M (Z',Z)})
\end{equation}
is given by twice  the integral of the  symplectic two-form $\omega_M$
of $M$  on the
surface  $\sigma (0,Z,Z')$ of the geodesic triangle
 $\gamma (0,Z,Z')\subset M$
 \begin{equation}\label{pha}
 \Phi_M (Z',Z)=2\int_{\sigma (0,Z,Z')}\omega_{M}.
 \end{equation}

Also
\begin{equation}\label{dcc}
|(\underline{{e}}_{Z'},\underline{{e}}_{Z})|=
|(\underline{{e}}_{i (Z')},\underline{{e}}_{i (Z)})|=
\cos d_C(i (Z'), i (Z)).
\end{equation}}
\end{theorem}
In the last relation
$$  d_C([u],[v])=\arccos\frac {|\skp {u}{v}|}{\Vert u\Vert\Vert v\Vert}\ .$$

Note that Theorem \ref{sb2}, which is Theorem 2b) and Theorem 3 in
\cite{sb7}, was obtained in \cite{clerc} under the form of Theorem 2.1.

We are interested to find other manifolds for which eq. (\ref{basic}) still
holds. But first we present some calculation of the
  the phase which appears when we multiply
two representations  on the complex Grassmann manifold.

%\section{ The two-cocycle on $SU(n,m)$ and $SU(n+m)$}\label{two}
\section{ Multiplicative factors on the complex Grassmann
manifold and its noncompact dual}\label{two}
\begin{proposition}\label{eu}
Let the noncompact (compact) Grassmann manifold
$X_{n}$ ($X_c$) parameter\-ized in the $B$ and $Z$ parametrizations
(\ref{57.1})-(\ref{57}), where the
 parameters are related
by equation (\ref{58}). Let $\sigma$ be  the section
which associates to the element in $Z\in\got{p}_+$ the element in $G_{n,c}$
given by the equation (\ref{79}). Then,  to a $Z\in\got{p}_+$
there corresponds a  $g\in G$ such that $g\cdot o=Z$. Let $D(B)$ represents
the matrix (\ref{57.1}) expressed in the form  (\ref{79}). Then
\begin{equation}\label{coc}
D(B_1)D(B_2)=D(B_3)\times e^{i\Phi}, 
\end{equation}
where the multiplicative factor $\Phi$  has the expression
\begin{equation}\label{faza}
  e^{i\Phi}=\det\left[ (\oo-\epsilon Z_1Z^+_2)(\oo-\epsilon
Z_2Z^+_1)^{-1}\right]^{-\epsilon/2}.%=e^{-2\pi i f(g_1,g_2)}
\end{equation}
%and the 2-cocycle $f$  has the expression (\ref{fcociclu}).
Here $Z_3=Z(B_3)$, where $Z_3$ is given by eq. (\ref{z3}).

\end{proposition}
{\it Sketch of the Proof}.

We give two proofs.

A) First proof  is a matrix calculation. We
 indicate the main steps.

a) The product $D(B_1)D(B_2)$ in the 
$Z_1,Z_2$-parametrization of the form (\ref{79}) is written down as
a four-block matrice 
\begin{equation}\label{mnpq}
D(B_1)D(B_2)=\left(\begin{matrix} M&N\cr P& Q\cr
\end{matrix}\right)
\end{equation}
where
\begin{subequations}
\begin{eqnarray}
M & = & (\oo +\epsilon Z_1Z_1^+)^{-1/2}(\oo - \epsilon Z_1Z_2^+)
(\oo +\epsilon Z_2 Z^+_2 )^{-1/2},\label{mmm} \\
N & = & (\oo +\epsilon Z_1Z^+_1)^{-1/2}(Z_2+Z_1)(\oo +\epsilon
Z_2^+Z_2)^{-1/2}, \label{nnn}\\
P & = &-\epsilon (\oo + \epsilon Z^+_1Z_1)^{-1/2}(Z_1+Z_2)^+
(\oo + \epsilon Z_2Z_2^+)^{-1/2},\label{ppp} \\
Q & = & (\oo + \epsilon Z^+_1Z_1)^{-1/2}(\oo -\epsilon Z^+_1Z_2)
(\oo +\epsilon Z_2^+Z_2)^{-1/2}.\label{qqq}
\end{eqnarray}
\end{subequations}
Here we present the calculation of $M$, the calculations of the other matrices
being similar:
\begin{equation*}
\begin{split}
M & =  (\oo +\epsilon Z_1Z_1^+)^{-1/2}(\oo + \epsilon Z_2Z^+_2)^{-1/2}-
\epsilon Z_1(\oo + \epsilon Z^+_1Z_1)^{-1/2}
(\oo + \epsilon Z^+_2Z_2)^{-1/2}Z^+_2 \\
 & =  (\oo +\epsilon Z_1Z_1^+)^{-1/2}(\oo + \epsilon Z_2Z^+_2)^{-1/2}-
\epsilon (\oo + \epsilon Z_1Z^+_1)^{-1/2}Z_1Z^+_2
(\oo + \epsilon Z_2Z^+_2)^{-1/2} \\
& =  (\oo +\epsilon Z_1Z_1^+)^{-1/2}(\oo - \epsilon Z_1Z_2^+)
(\oo +\epsilon Z_2 Z^+_2 )^{-1/2}. 
\end{split}
\end{equation*}
b) We write down again a Gauss decomposition of the product in eq. (\ref{mnpq})
\begin{equation}\label{gauss}
\left(\begin{matrix}M & N\cr P & Q\cr\end{matrix}\right)=
\left(\begin{matrix}\oo & Z'\cr {\bf 0} &\oo\cr \end{matrix}\right)
\left(\begin{matrix}\alpha & \bf{0}\cr \bf{0}&  \beta\cr
\end{matrix}\right) \left(\begin{matrix}\oo & \bf{0}\cr Z& \oo\cr
 \end{matrix}\right),
\end{equation}
where 
\begin{subequations}
\begin{eqnarray}%\label{MNPQ}
M  & = &  \alpha + Z'\beta Z,\label{MM}\\
N  &  = &   Z'\beta,\label{NN}\\ 
P  & =  & \beta Z,\label{PP}\\
Q  &  =  & \beta .\label{QQ}
\end{eqnarray}
\end{subequations}
It results $Z'=NQ^{-1}$, and finally, and it is find that 
  $Z'\equiv Z_3$, where
$Z_3$ has the expression given by  eq. (\ref{z3}). 
Now we calculate $\alpha$ in the Gauss decomposition  (\ref{gauss}), using
eq. (\ref{MM}): 
\begin{equation}\label{mare}
\begin{split}
\alpha & = M - NQ^{-1}P =  
 (\oo +\epsilon Z_1Z_1^+)^{-1/2}(\oo - \epsilon Z_1Z_2^+)
(\oo +\epsilon Z_2 Z^+_2 )^{-1/2}+\\
& \epsilon (\oo+\epsilon Z_1Z_1^+)^{-1/2}
(Z_1+Z_2)(\oo+ \epsilon Z^+_2Z_2)^{-1/2}(\oo + \epsilon Z^+_2Z_2)^{1/2}
(\oo - \epsilon Z^+_1Z_2)^{-1}\times\\
 & (\oo + \epsilon Z^+_1Z_1)^{1/2}(\oo + \epsilon Z^+_1 Z_1)^{-1/2}
(Z_1+Z_2)^+(\oo + \epsilon Z_2Z_2^+)^{-1/2},
\end{split}
\end{equation}
so we have for $\alpha$
\begin{equation}\label{alfa}
\alpha = (\oo + \epsilon Z_1 Z_1^+)^{-1/2}\Lambda (\oo + 
\epsilon Z_2Z_2^+)^{-1/2},
\end{equation}
where
\begin{equation}\label{lala}
\Lambda = \oo -\epsilon Z_1Z^+_2 +
\epsilon (Z_1+Z_2)(\oo - \epsilon Z^+_1Z_2)^{-1}(Z_1+Z_2)^+.
\end{equation}
In order to find a simpler expression for $\Lambda$,
 we substitute in eq. (\ref{lala}) for $Z_1+Z_2$
\begin{subequations}
\begin{align}
Z_1+ Z_2 &= Z_1(\oo - \epsilon Z_1^+Z_2)+(\oo +\epsilon Z_1Z_1^+)Z_2, \\
\intertext{and, similarly for its adjoint}
(Z_1+Z_2)^+ & = (\oo - \epsilon Z^+_1Z_2)Z^+_2+Z^+_1(\oo + \epsilon Z_2Z_2^+).
\end{align}
\end{subequations}
For $\Lambda$ in eq. (\ref{alfa}) we finally find
\begin{equation}\label{lambda}
\Lambda = (\oo +\epsilon Z_1 Z_1^+)(\oo -\epsilon Z_2 Z_1^+)^{-1}
(\oo +\epsilon Z_2 Z_2^+).
\end{equation}
c) In accord  with the Gauss decomposition (\ref{56}) of equation (\ref{79}), 
we have 
$$ D(B_3)= \left(\begin{matrix} \oo & Z_3\cr \bf{0}& \oo
\end{matrix} \right) \left(\begin{matrix}
U & \bf{0}\cr \bf{0} & V\cr \end{matrix} \right)
 \left(\begin{matrix}\oo & \bf{0}\cr 
-\epsilon Z^+_3 & \oo\end{matrix} \right),       $$
\begin{equation}\label{U}
U  =  (\oo + \epsilon Z_3Z_3^+)^{1/2},
\end{equation}
where $Z_3$ has the expression given by eq. (\ref{z3}). The calculation
is similar with equations  (8.6) and (8.7)  in \cite{grass},
  where in
eq. (8.3) we have to substitute $Z_1\rightarrow -Z_1$. This corresponds
to the  fact:  $U^{-1}(Z)=
U(-Z)$ (see the Appendix).

The final expression  needed to determine $U$ in eq. (\ref{U}) is:
\begin{eqnarray}
\nonumber
 \oo +\epsilon Z_3Z_3^+ & = & (\oo +\epsilon Z_1Z_1^+)^{1/2}
(\oo-\epsilon Z_2 Z_1^+)^{-1}\times\\
\label{zzz}
 & & (\oo +\epsilon Z_2Z_2^+)(\oo -\epsilon Z_1Z_2^+)^{-1}
(\oo + \epsilon Z_1 Z_1^+)^{1/2}.
\end{eqnarray}

d) We use the following  relation (see the case of the
maximal weight for the compact case  e.g. eq. (3.12)
in \cite{bg}) for the action of the representation ${\pi}$ on the
minimal  weight vector $o$:
\begin{equation}\label{pondere}
j_1=j_2 =\cdots j_n=k,j_{n+1}= j_{n+2}=\cdots
j_{n+m}=0, 
\end{equation}
\begin{equation}\label{asaofi}
 {\pi}\left( \begin{matrix}A &  0\cr   0 & D \cr
\end{matrix}\right) o = (\det A)^{-k\epsilon} ~o,~ \det A\det D =1  .
\end{equation}
The sphere $SU(2)/U(1)$ corresponds to $k=2j, j=\frac{1}{2}, 1,\cdots$,
while the dual noncompact case $SU(1,1)/U(1)$ is obtained replacing
$k\rightarrow 2k=2,3,\cdots$. For the case of the complex
Grassmann manifold and its noncompact dual below we take $k=1$ if not
specified otherwise.
The phase $\Phi$ in eq. (\ref{coc}) is obtained from the relation
$$(\det \alpha)^{-\epsilon} = (\det U)^{-\epsilon} e^{i\Phi},$$
where $\alpha$ is given by the equations (\ref{alfa}), (\ref{lambda})
and $U$ by equations (\ref{U}), (\ref{zzz}). ~ $\gata$

B). We now present briefly a second proof of Proposition \ref{eu}.
  We recall firstly some general considerations on multipliers and
coherent states \cite{tim}, \cite{buc}. Here again $\pi$ is unitary
irreducible
representation of the group $G$ on a Hilbert space \Hi .

%Recalling the  definition (\ref{rel1})  (respectively, (\ref{sect2})
%of the function $f'_{\psi}$  
%($f_{\psi}$), we have
Let 
\begin{equation}\label{iar}
f_{\psi}(z)=(e_{\bar{z}},\psi)=\frac{(\pi(\bar{g})e_0,\psi)}
{(\pi(\bar{g})e_0,e_0)},~z\in M, \psi\in\Hi .
\end{equation}
We get 
\begin{equation}\label{iar1}
 f_{\pi(\overline{g'}).\psi}(z)= \mu
(g',z)f_{\psi}(\overline{g'}^{-1}.z),
\end{equation}
where 
\begin{equation}\label{iar2}
 \mu (g',z)=
\frac{(\pi(\overline{g'}^{-1}\overline{g})e_0,e_0)}{(\pi(\overline{g})e_0,e_0)}
=\frac{\Lambda (g'^{-1}g)}{\Lambda (g)}.
\end{equation}
We recall that 
\begin{equation}\label{FazaP}
{\pi(g).e_0=e^{i\alpha (g)}e_{\tilde{g}}
=\Lambda (g) e_{z_{g}}}\end{equation}
where we have used the decompositions,
$g=\tilde{g}.h, ~(G=G/H.H);~~ g = z_g. b ~(G_{\C}=G_{\C}/B.B)$.
 We have also  the relation
$\chi_0(h)=e^{i\alpha (h)},~ h\in H$ and 
$ \chi (b) =\Lambda (b),~b\in B$, where
$\Lambda(g)=\frac{e^{i\alpha(g)}}{(e_{\bar{z}},e_{\bar{z}})^{1/2}}$.
We can also write down another expression for the multiplicative
factor $\mu$ appearing in  eq. (\ref{iar1}) using CS-vectors 
\begin{equation}\label{iar3}
\mu (g',z)=\Lambda(\bar{g'})(e_{\bar{z}},e_{\bar{z'}})=e^{i\alpha
  (\bar{g'})}
\frac{(e_{\bar{z}},e_{\bar{z'}})}{(e_{\bar{z'}},e_{\bar{z'}})^{1/2}}
\end{equation}
and $$\arg\mu (g',z)= \alpha
  (\bar{g'})+ \Phi_M(\bar{z},\bar{z'}). $$
The following assertion is easy to be checked using successively
eq. (\ref{iar2}):
\begin{Remark}
\label{unau}Let us consider the relation (\ref{iar}). Then we have
(\ref{iar1}), where  $\mu$ can be written down as in equations
(\ref{iar2}), 
(\ref{iar3}). 
We have the relation
\begin{equation}\label{incauna}
 \mu (g,z) =J(g^{-1},z)^{-1},
\end{equation} i.e. the multiplier $\mu$
 is the cocycle in the  unitary representation
$(\pi_K,\Hi_K)$ attached to the
positive definite holomorphic kernel 
$K(z,\bar{w}):=(e_{\bar{z}},e_{\bar{w}})$ % defined by equation (\ref{kernel}),
 \begin{equation}\label{num}
(\pi_K(g).f)(x):=J(g^{-1},x)^{-1}.f(g^{-1}.x),
\end{equation}
and the cocycle verifies the relation
\begin{equation}\label{prod}
J(g_1g_2,z)=J(g_1,g_2z)J(g_2,z).
\end{equation}
\end{Remark}
  Note that
{\em   the prescription (\ref{num})
 defines
a continuous action of $G$ on ${\mbox{\rm{Hol}}}(M,\C )$ with respect to
the compact open topology on the space  ${\mbox{\rm{Hol}}}(M,\C )$.
 If $K:M\times \bar{M}\rightarrow \C$ is a continuous positive definite kernel
holomorphic in the first argument satisfying
\begin{equation}\label{relat}
K(g.x,\overline{g.y})=J(g,x)K(x,\bar{y})J(g,y)^*,
\end{equation}
$g\in G$, $x,y\in M$, then the action of $G$ leaves the reproducing kernel
Hilbert space $\Hi_K\subseteq {\mbox{\rm{Hol}}}(M,\C )$ invariant and defines
a continuous unitary representation $(\pi_K,\Hi_K)$ on this space}
(cf. Prop. IV.1.9 p. 104 in  Ref. \cite{neeb}).

In Perelomov's notation at p.42 in \cite{per1}, eqs. (\ref{iar1})-(\ref{iar3})
which define the multiplicative factor $\mu$ read
 \begin{equation}\label{per}
{\pi}(g_1)e_B= e^{i\beta (g_1,z)}e_{g_1.B}.
\end{equation}
Let us recall the notation $(e_z,e_{z'})=e^{i\Phi
(z,z')}|(e_z,e_{z'})|$.
With the
Remark \ref{unau}, it is easy to see that 
\begin{equation}\label{dif}
\beta (g_1,z)=\Phi (z_{g^{-1}_1},z_g)-\alpha (g^{-1}_1).
\end{equation}
 It can be proved that the phase $\alpha $ on
$X_{c,n}=G_{c,n}/S(U(n)\times U(m)$ is given by
\begin{equation}\label{incuna}
{\pi}\left(\begin{matrix} A_1 & B_1\cr  C_1 & D_1
 \cr\end{matrix}\right)e_0=
\left[\frac{\det (A_1)}{\overline{\det}
 (A_1)}\right]^{-k\frac{\epsilon}{2} }e_{B},
~\left(\begin{matrix} A_1 & B_1\cr  C_1 & D_1
 \cr\end{matrix}\right)\in G_{c,n}.
\end{equation}
$B\in\got{g}_{c,n}$ in eq. \ref{incuna}
 is given by eq. (\ref{invers}) with $Z=B_1D^{-1}_1$ or $$ZZ^{+}
= \epsilon [(A_1A_1^+)^{-1}-1] . $$
But if the matrix of the group $G_c$ ($G_n$) is taken from the
 homogeneous space $X_c$ (respectively, $X_n$),
then in eq. (\ref{dif}) $\alpha =0$. 
Then it is used the relation $(e_z,e_{z'})=\det (\oo +
z'z^+)^{\epsilon}$ and eq. (\ref{faza}) is re-obtained.

It can be seen that eq. (\ref{faza}) can be deduced from the equation:
\begin{equation}
\label{marea}
{\pi}\left(\begin{matrix} A_1 & B_1\cr C_1 & D_1\cr\end{matrix}\right)e_B=
\left[\frac{\det (A_1-\epsilon B_1Z^+)}{\overline{\det} (A_1-\epsilon
B_1Z^+)}\right]^{-k\frac{\epsilon}{2} }e_{g_1.B}.
\end{equation}
We take $k=1$ and the matrix $g_1$ is taken of the form given by eq.
 (\ref{79}). $\gata$

So, in this section we have presented a brute-force calculation
 of the multiplicative factor (\ref{coc}) and a simpler proof of the
same
calculation. In the next section the
 exact meaning of this cocycle will be clarified.

We end the section with another 

\begin{Remark}\label{unuc}
The relation (\ref{relat}) can be used to find the cocycle $J$ in
equation (\ref{num}).
\end{Remark}

We illustrate this assertion by the
 example of  the complex Grassmann manifold \newline \Gras~   and its
noncompact dual.
  The scalar product (the reproducing kernel) corresponding to the
extremal weight (\ref{pondere})  is 
\begin{equation}\label{scgras}
K(X,\bar{Y})=(e_{\bar{X}},e_{\bar{Y}})=\det (\oo +\epsilon XY^+)^{\epsilon k}.
\end{equation}
Below we take $k=1$. Then
$$
K(g.X,g.Y)=\det (XB^+-\epsilon A^+)^{-\epsilon}K(X,Y)\overline{\det}(YB^+-
\epsilon A^+)^{-\epsilon } .$$
$$J(\left(\begin{matrix} A&  B\cr
C & D\cr \end{matrix}\right),X)=\det (A^+-\epsilon XB^+)^{-\epsilon}$$
is an automorphy factor (see e.g.  \cite{cart}) and eq. (\ref{prod}) is
satisfied. 
If $$g=\left(\begin{matrix} A&  B\cr
C & D\cr \end{matrix}\right)\in SU(n+m)~ \text{or}~ SU(n,m) ,$$
then $$g^{-1}=\left(\begin{matrix} A^+&  \epsilon C^+\cr
\epsilon B^+ & D^+\cr \end{matrix}\right) ,$$
and 
\begin{equation}\label{auto}
J(\left(\begin{matrix} A&  B\cr
C & D\cr \end{matrix}\right)^{-1},X)=\det(A-XC)^{-\epsilon } .
\end{equation}
 Equation (\ref{num}) reads in this case
$$
{\pi}\left(\begin{matrix} A&  B\cr
C & D\cr \end{matrix}\right)f(X)=\det (A-XC)^{\epsilon}f[(A^+X +
\epsilon C^+)(\epsilon B^+X+D^+)^{-1}].
~\gata $$
\section{Two-cocycles and symplectic areas of geodesic \\ triangles
 on hermitian symmetric spaces}\label{compar}

Firstly we review some results obtained by A. Guichardet and
D.  Wigner, and  J.-L. Dupont and
A.  Guichardet. Then, using their results, we answer partially to the question
addressed at the end of  \S  \ref{question}. 

a)
A. Guichardet and D. Wigner \cite{guit}
 have  proved that: {\em  a simple Lie group has non-trivial
 continuous 2-cohomology group $H^2(G, \R )$ if and only if $G/K$ admits a
 $G$-invariant complex structure},  where $K$ is  a maximal compact subgroup of
 $G$.
More exactly, the starting point of their investigation is the following 
lemma (presented here in abbreviated form), based mostly 
on the results collected in Helgason's book \cite{helg}:
\begin{lemma}
Let $G$ be simple. Then:\\
(a) $[\got{p},\got{p}]=\got{k}$;\\
(b) the adjoint representation of $\got{k}$ in $\got{p}$ is irreducible;\\
(c) the next conditions are equivalent:\\
(i) $H^2(G,\R )\not= 0$;\\
(ii) ${\mbox{\rm{Hom}}}_{\got{k}}\,(\bigwedge^2\got{p},\R )\not= 0$;\\
(iii) there is a $\got{k}$-invariant complex structure;\\
(iv) $G/K$ admits a  $G$-invariant complex structure;\\
(v) the center $\got{Z}(\got{k})$ of $\got{k}$ is non void;\\
(vi) ${\mbox{\rm{Hom}}}\,(\got{k},\R )\not= 0$;\\
(d) If the previous conditions are fulfilled, then:\\
(i') ${\mbox{\rm{dim}}}\, H^2(G,  \R )= {\mbox{\rm{dim}}}\, \got{Z}(\got{k})
=  {\mbox{\rm{dim}}}\, {\mbox{\rm{Hom}}}\,(\got{k},\R)=1$;\\
(ii') the $G$-invariant complex structure on $G/K$ is hermitian.
\end{lemma}
Above  $\got{g}= \got{k}+\got{p}$ is a Cartan decomposition. 

Here is the well known list of groups $G$ which have real
hermitian Lie algebras: $SU(m,n)$; $SO_0(2,q)$; ($q=1$ or
$q\ge 3$); $Sp(n,\R )$ $(n\ge 1)$; $SO^*(2n)$ ($n\ge 2$); $E_6$; $E_7$. 

A. Guichardet and D.  Wigner have considered the  real differentiable
 2-cocycle $f$:
 \begin{equation}
\label{2coc}
 f(g_1,g_2)=\frac{1}{2\pi}{\mbox{\rm{arg}}} (v(g_1)v(g_2)v(g_1g_2)^{-1}),
  \end{equation}
where  $v $ is a non-trivial  morphism   of $G$ in the torus $T$.

Note that $H^s(G,\R ) = H^s(\got{g},\got{k},\R )= H^s(\hat{G}/K , \R )=
\text {Hom}_K(\bigwedge^s\got{p},\R )$,
where $\hat{G}$ is the compact form of the real noncompact simple
 Lie group $G$ (cf. \cite{guit1}). More exactly, 
cf. Proposition 7.6 in \cite{guit1}:
 {\em $f(g_1,g_2)$ defines a two-cocycle $f\in Z^2_{diff}(G,\R )_K$
whose class in $H^2(G,\C)$ corresponds to the element in 
$H^2(\got{g},\got{k},\C )$ via the van Est isomorphism.}
\begin{Remark}\label{ident1}
The multiplicative factor $\Phi$ in eq. (\ref{faza}) is the
two-cocycle (\ref{2coc}) determined by A. Guichardet and 
D. Wigner for the group $G= G_n= SU(n,m)$
expressed in Pontrjagin's coordinates on $X_n$ ($\epsilon = -1$).
\end{Remark}
{\em Proof}.
For  
  $g\in G_n$ put in the block form
 \begin{equation}
\label{block}
 g= \left(\begin{array}{cc}a & b\\  c & d\end{array}\right)
  \end{equation}
it  has been  shown in \cite{guit}
  that $v$ in eq. (\ref{2coc}) is given by  $v(g)=\det a$.
In  Pontrjagin's coordinates $Z=\pi (g)$, ($\pi\sigma =1$, where
$\pi$ is the natural projection $G_{c,n}\rightarrow X_{c,n}$),
the function $v$ in   the cocycle (\ref{2coc}) is (below, in the
formulas
from \cite{guit} $\epsilon =-1$):
 \begin{equation}
 v(Z)= v(\pi (g))= \det ({\oo} + \epsilon ZZ^+)^{-1/2}.\label{vv}
 \end{equation}

     Multiplying two matrices $g_i=g(Z_i)$, $i=1,2$
 of the type (\ref{block}) we get a matrix of the
 same form, where the block of the type $a$ is the $M$ given by eq. 
(\ref{mmm}), i.e. 
 \begin{equation}\label{mult}
a= (\oo + \epsilon Z_1Z^+_1)^{-1/2}(\oo - \epsilon Z_1 Z^+_2)(\oo + 
 \epsilon Z_2 Z^+_2)^{-1/2}.
 \end{equation}
Combining eqs. (\ref{vv}) and (\ref{mult}), we get for the
 2-cocycle (\ref{2coc}) the expression 
 \begin{equation}
f(Z_1,Z_2)=\frac{1}{2\pi}{\mbox{\rm{arg}}}\det (\oo - \epsilon Z_1Z^+_2)^{-1},
  \end{equation}
 \begin{equation}\label{fcociclu}
 f(Z_1,Z_2)=\frac{1}{4\pi i}\log \frac {\det (\oo - \epsilon Z_2Z^+_1)}
 {\det (\oo - \epsilon Z_1 Z^+_2)}. ~~~\gata
\end{equation} 
So, Proposition \ref{eu} gives the expression (\ref{fcociclu}),
independent of the results of A. Guichardet and D. Wigner, in Pontrjagin's
coordinates. More exactly, 
\begin{equation}\label{gutu}
 e^{i\epsilon \Phi}=e^{-2\pi i f(g_1,g_2)}.
\end{equation}

b)
The geometric significance of the 2-cocycle (\ref{2coc}) was found by
J.-L.  Dupont and A. Guichardet \cite{dupont}. Let in their notation
 $v_{*e}$ be the differential of the
 homomorphism $v$ at the origin $e\in G$ and $P=v_{*e}/2\pi i$, and let
 $P(\Omega )$ be the $G-$invariant differential 2-form on $G/K$ with the
 value
 of the origin $o$   given by $P(\Omega )_0(A,B)=-\frac{1}{2}P([A,B])$, 
 $A,B\in \got{p}_n$. Let 
$\Delta (g_1, g_2)$ be the
 geodesic cone with corner $o$ and base the geodesic joining $g_1\cdot o$ with
 $g_1
 g_2\cdot  o$. Then, in a previous publication
 (see references in \cite{dupont})
J.-L. Dupont
 has constructed the 2-cocycles by
integration of $G-$invariant differential forms on geodesic simplexes
in a symmetric space $G/K$, where $K$ is a maximal compact subgroup of
$G$
 \begin{equation}
\label{dg}
 c(g_1,g_2)=\int_{\Delta (g_1,g_2)}P(\Omega ).
 \end{equation}

In the quoted paper \cite{dupont} J.-L.  Dupont and A. Guichardet have
proved the equality 
\begin{theorem}\label{egal}
\begin{equation}
f(g_1,g_2)=c(g_1,g_2),~ g_1,g_2\in G.
\end{equation}
\end{theorem}

\begin{Remark}\label{ident2}Now we only express 
  Theorem \ref{egal} in Pontrjagin's coordinates.
\end{Remark}
 In order to calculate $c(g_1,g_2)$, we have to calculate
$\got{I}(0,Z_1,Z_3)$ with formula (\ref{ar}), i.e.
\begin{equation}
c(g_1,g_2)=\frac{\epsilon}{4 i}\log A{(A^+)^{-1}},
\end{equation}
where
$$A= \oo +\epsilon Z_1Z_3^+,$$
and $Z_3$ is given by eq. (\ref{z3}).
It is obtained
$$A=(\oo + \epsilon Z_1Z_1^+)^{1/2}B
( \oo + \epsilon Z_1Z_1^+)^{-1/2},$$ while for $B$
it is obtained
$$B = (\oo - \epsilon Z_1Z_2^+)^{-1}(\oo +\epsilon Z_1Z_1^+).$$
We get 
\begin{equation}
c(g_1,g_2)=\frac{\epsilon}{4 i}\log\det \frac{(\oo-\epsilon Z_2Z^+_1)}
{(\oo -\epsilon
Z_1Z_2^+)}.
\end{equation}
We checked the validity of Theorem \ref{egal} and get
 $f(g_1,g_2)=
\frac{\epsilon}{\pi}c(g_1,g_2)$, the multiplicative factor coming
 from  a  different normalization. $\gata$

In Theorem \ref{gr} we have calculated the expression of the area
of the geodesic triangle with vertices $(0, z, z')$, while eq. 
(\ref{fcociclu}) gives the expression of the 2-cocycle $f$.

Putting together the results of the papers of A. Guichardet and D. Wigner,
and J.-L. Dupont and A. Guichardet \cite{guit,dupont} (see
 also the book of A. Guichardet
\cite{guit1}),
 we get an answer to our question of generalizing the relation (\ref{basic})
 to simple Lie groups:

\begin{theorem}\label{result}
Let $G $ be a simple Lie group and $K$ a maximal compact subgroup. 
Then the only  coherent state manifolds
$G/K$ for which the phase of the scalar product of two
coherent state vectors is twice the symplectic area 
of a geodesic triangle are the hermitian symmetric spaces.
\end{theorem}

So far, we have seen that between the {\bf simple Lie groups}
$G$,  {\bf only}
 the groups which have a real  hermitian simple Lie algebra
lead to coherent states based on $G/K$ which
 have the property (\ref{basic}). Meanwhile, the same property 
 (\ref{basic}) is verified by the Heisenberg group,
 as was underlined in \S \ref{br}.
Hence, it is natural to formulate the following 
{\bf question: For which groups $G$ which   admits coherent state
representations  the property (\ref{basic})
is still  true?}
 We recall that the group $G$ admits 
coherent state representations
 (cf. W. Lisiecki \cite{lis1} and  K. Neeb \cite{neeb}) if $G/H, ~H\subset K$ 
admits a holomorphic embedding in a projective
Hilbert space, where   $H$ is isotropy subgroup of the representation
with extreme weight vector 
$e_0$. For example,
property (\ref{basic}) is still true if
$G$    is a semi-direct product of a hermitian 
 type group (i.e. $G/K$ is hermitian symmetric) and a Heisenberg group? The
answer to this question is given by those CS-groups $G$
which lead to naturally reductive homogeneous spaces $G/H$.

In this context of Theorem \ref{result},
 we would like to recall our result established in \cite{tim,buc}:

\begin{Remark}\label{red}  The coherent state manifold $M\cong G/H$, for
which the isotropy representation has discrete kernel, or for
admissible Lie algebras and faithful CS-representations, is a
reductive space.
\end{Remark}

In the same context, let us recall the following classical result
\cite{hano}:
\begin{theorem}
Let $G/B$ be a K\"ahlerian homogeneous space of a reductive Lie group
$G$ and let $G$ be effective on $G/B$. If the Riemannian connection on
$G/B$ induced by the invariant K\"ahlerian metric is the canonical
affine connection of the first kind with respect to a certain
$B$-invariant decomposition of $\got{g}$, then $G/B$ is hermitian symmetric.
\end{theorem}

\setcounter{equation}{0}\renewcommand{\theequation}{A.\arabic{equation}}
\section{Appendix: parametrization of
the Grassmann ma\-nifold and its noncompact dual}\label{app}

The elements $U\in G_{n,c}$ verify the relation
\begin{equation}\label{rule}
U^+I_{nm}(\epsilon )U=I_{nm}(\epsilon ),~~I_{nm}(\epsilon )=
\left(\begin{matrix} \epsilon \oo_n & \bf{0} \cr
\bf{O} & \oo_m \cr \end{matrix}\right), 
\end{equation}
where $\epsilon =1 (-1)$ corresponds to $G_c$ (resp. $G_n$).

Also we have the Cartan decomposition
\begin{equation}\label{cartan}
\got{g}_{n,c}=\got{k}+\got{p}_{n,c}; 
\end{equation}
$\got{g}_{n,c}$ ($\got{k}$)
 denotes the Lie algebra of the group $G_{n,c}$
(respectively $K$), and
we have, locally (globally),
 the diffeomorphism of $X_{c}$ with $\got{p}_{c}$ (respectively, 
 $X_{n}$) with $\got{p}_{n}$
$$X_{n,c}=\exp({\got{p}_{n,c}}) o.$$

The manifold $X_c$  and its noncompact dual $X_n$ is parametrized by
$B\in\got{p}_{n,c}$
\begin{subequations}
\begin{eqnarray}
X_{n,c} & = & \exp\left(
\begin{matrix}{\bf 0}&B\cr
                           -\epsilon B^+&{\bf 0} \cr 
\end{matrix} \right) o=
   \left(\begin{matrix}{\rm co}\sqrt{BB^+}&B{\displaystyle 
\frac{{\rm si} \sqrt{B^+B}}{ \sqrt{B^+B}}}\cr
                            -\epsilon{\displaystyle \frac
{{\rm si} \sqrt{B^+B}}{
   \sqrt{B^+B}}}B^+&{\rm co}\sqrt{B^+B}\cr\end{matrix}\right)o
\label{57.1} \\*[2.0ex]
 \nonumber &  & \\
  & = &
 \left(\begin{matrix}{\oo}&Z\cr {\bf 0}&{\oo}\cr\end{matrix}\right)
\left(\begin{matrix}({\oo}+\epsilon ZZ^+)^{1/2}&
{\bf 0}\cr {\bf 0}&({\oo}+\epsilon Z^+Z)^{-1/2}\cr\end{matrix}\right)
\left(\begin{matrix}{\oo}&{\bf 0}\cr -\epsilon Z^+&{\oo}\cr\end{matrix}
\right)o
\label{56}\\*[1.0ex]
 \nonumber &  & \\
  & = & \left(\begin{matrix}(\oo + \epsilon ZZ^+)^{-1/2} &
 Z(\oo +\epsilon Z^+Z)^{-1/2}\cr -\epsilon (\oo + \epsilon Z^+Z)^{-1/2}
 Z^+ & (\oo +\epsilon Z^+ Z)^{-1/2} \cr \end{matrix}\right) o,
 \label{79} \\[1.0ex]
\nonumber &  & \\
 & = & \exp\left(\begin{matrix}{\bf 0}&Z\cr {\bf 0}&{\bf 0}\cr\end{matrix}
\right)P ,
\label{57}
\end{eqnarray}
\end{subequations}
where $\epsilon = 1 (-1)$ for compact (respectively non-compact) manifolds.
Here  {\rm co}  is the  circular cosine  cos 
(resp. the hyperbolic cosine coh) for $X_c$ (resp. $X_n$) and
 similarly for
 si.
 $Z$ is the  $n\times m$ matrix of Pontrjagin's
coordinates in ${\cal V}_0$ related to $B$ by the formula
\begin{equation}
Z=Z(B)=B{\frac{{\rm ta}\sqrt{B^+B}}{ \sqrt{B^+B}}}~ ,\label{58}
\end{equation}
and ta - the hyperbolic tangent tanh (resp. the
circular tangent tan) for $X_n$ (resp. $X_c$). The relation inverse to
eq. (\ref{58}) is
\begin{equation}\label{invers}
B=\frac{\text{arcta}\sqrt {ZZ^+}}{\sqrt {ZZ^+}}
\end{equation}

The transitive action of an element of the group $G_c=SU(n+m)$ 
($ G_n=SU(n,m)$) on $X_c$ (resp. $X_n$) is given by the linear
fractional transformation
\begin{equation}\label{fractionar}
Z'=Z'(Z)=U\cdot Z = (AZ+B)(CZ+D)^{-1},~ U =\left(\begin{array}{cc}
A & B\\ C & D \end{array}\right)\in G_c(G_n).
\end{equation}

So, we have to find a matrix $U\in G_c(G_n)$ such that eq. (\ref{rule})
is satisfied, i.e.
\begin{equation}\label{prod1}
\left\{
\begin{array}{ccc}
A^+A+\epsilon C^+C & = & \oo_n,\\
\epsilon B^+ B + D^+D & = & \oo_m, \\
\epsilon B^+ A + D^+ C & = & \bf{0}.
\end{array}
\right.
\end{equation}

Now let $g_i\cdot 0=Z_i$, $i=1,2$, and let 
\begin{equation}\label{action}
Z_3:=g_1\cdot Z_2=\left(\begin{matrix} A_1 & B_1\cr C_1 & D_1\cr
 \end{matrix}\right)(Z_1)\cdot Z_2=
(A_1Z_2+B_1)(C_1Z_2+D_1)^{-1}.
\end{equation}
The matrix appearing in eq. (\ref{79}),
a solution of  eq. (\ref{prod1}), fixes  a section $\sigma:
G/K\rightarrow G$ such that $\sigma (o)= e$.  It can 
be shown (similar to eq. (8.3) in \cite{grass},
but with $Z_1\rightarrow -Z_1$) that $Z_3$ has the expression 
(\ref{z3}):
 \begin{equation}\label{z3}
Z_3=(\oo +\epsilon Z_1Z_1^+)^{-1/2}(Z_1+Z_2)(\oo -\epsilon Z_1^+Z_2)^{-1}
(\oo + \epsilon Z_1^+Z_1)^{1/2}.
\end{equation}

We mention also the following useful relation 
(which enables to get eq. (\ref{z3})
 using eq. (8.3) in \cite{grass}):
{\em Let $\sigma : X_{n,c}
\rightarrow G_{n,c}$ be the section with the property
 that $\sigma (o)=e$ which 
associates to a point in the Grassmann manifold $X_{n,c}$ the matrix $U$
  (\ref{fractionar})  given by (\ref{79}) such that $U\cdot o = Z\in X_{n,c}$. 
Then $U^{-1}(Z)= U(-Z)$.}\\[3ex]

{\bf{ Acknowledgments}} I would like to thank the organizers of the {\it XX
Workshop on Geometric Methods in Physics} in  Bia\l owie\.{z}a, Poland 
 and to the organizers of the {\it Sixth International Workshop on differential
geometry and its applications} in Cluj, Romania
 for inviting me to
present  talks on this subject. Also I thanks Jiri Tolar and Alan Weinstein
for their comments and suggestions 
to my  talks at the seminaries 
  in their groups on the same theme.
 I am grateful  to Jean-Louis Clerc, Johan Dupont,  Bernard Magneron,  
Domingo Toledo and Mariano Santander for bringing in my attention their papers.
Discussions with Martin Schlichenmaier during my visit in Mannheim 
under the DFG-Romanian Academy Project 436 Rum 113/15  are acknowledged. 

\end{document}